\documentclass[11pt]{amsart}
\usepackage{amsfonts,amsmath, eucal, amssymb, upref}

\usepackage{chatterjee1}

\begin{document}
\title[Concentration of Haar measures]{Concentration of Haar measures, with an application to random matrices}
\author{Sourav Chatterjee}
\address{\newline367 Evans Hall \#3860\newline
Department of Statistics\newline
University of California\newline
Berkeley, CA 94720-3860\newline
{\it E-mail: \tt sourav@stat.berkeley.edu}\newline 
{\it URL: \tt http://www.stat.berkeley.edu/$\sim$sourav}
}
\subjclass[2000]{60E15, 28C10, 46L54, 15A52}
\keywords{Concentration of measure, concentration inequalities, Stein's method, semigroup method, Haar measure, random walk, mixing time}
\maketitle

\begin{abstract}
We show that the mixing times of random walks on compact groups can be used to obtain concentration inequalities for the respective Haar measures. As an application, we derive a concentration inequality for the empirical distribution of eigenvalues of sums of random hermitian matrices, with possible applications in free probability. The advantage over existing techniques is that the new method can deal with functions that are non-Lipschitz or even discontinuous with respect to the usual~metrics.
\end{abstract}

\section{Introduction and results}\label{intro}
Much attention has been paid to the derivation of concentration inequalities through logarithmic Sobolev inequalities, semigroup or transportation methods. Let us refer to the monograph of Ledoux \cite{ledoux01} for an extensive survey. On the other hand, starting with the pioneering work of Diaconis and Saloff-Coste \cite{ds96a, ds96} (see also \cite{saloffcoste97}), it is known that the rate of convergence to equilibrium of certain ergodic Markov semigroups or of random walks on compact groups involve logarithmic Sobolev constants. In this paper we make an explicit connection between the mixing time of a random walk on a compact group and the concentration property of the Haar measure. In other 
words the rate of convergence to equilibrium and the rate of concentration are directly connected. This is made precise in Theorem \ref{gpaction} below. 

We demonstrate that this new approach to concentration can 
be succesfully applied to random matrices. The gain as compared to previous methods is that it allows to deal with functions that are possibly non-Lipschitz  with respect to the Hilbert-Schmidt norm. 

The new method can be called an extension of Stein's method of exchangeable pairs \cite{stein86}, as developed by the author in his Ph.D.\ thesis \cite{chatterjee05}. From a different angle, it can also be viewed as a discrete analog of the semigroup tool for measure concentration (see Ledoux~\cite{ledoux01}, section on `semigroup tools'; see also the discussion in Subsection \ref{outline} of this paper). Some other applications of our method can be found in \cite{chatterjee06}.

The paper is organized as follows. We begin with the random matrix example (Theorem \ref{free}), followed by the statement of the main result (Theorem \ref{gpaction}), and a sketch of the proofs (Subsection \ref{outline}). Section \ref{intro} ends with a brief discussion of the literature. Proofs of Theorems \ref{gpaction} and \ref{free} are given in Section~\ref{proof}.

\subsection{A random matrix example} 
Let $M$ be an $n \times n$ complex hermitian (i.e.\ self-adjoint) matrix. The following terminology is standard.
\begin{enumerate}
\item The {\it empirical spectral measure} of $M$ is the probability measure on $\rr$, denoted by $\mu_M$, which puts $1/n$ on each eigenvalue of $M$, repeated by multiplicities.
\item The {\it empirical distribution function} of $M$,  denoted by $F_M$, is the cumulative probability distribution function corresponding to the empirical spectral measure.
\item Any hermitian matrix that has the same spectrum as $M$ can be written as $UM U^*$ for some unitary matrix $U$. Thus, the Haar measure on the group of unitary matrices of order $n$ naturally induces a `uniform distribution' on the set of all hermitian matrices with the same spectrum as $M$. We denote this probability measure by $\rho_M$. 
\end{enumerate}
\begin{thm}\label{free}
Let $M$ and $N$ be two hermitian matrices of order $n$. Suppose $A \sim \rho_M$ and $B\sim \rho_N$ are two independent random hermitian matrices. Let $H = A+B$. Then, for every $x\in \rr$,
$\var(F_H(x)) \le \kappa n^{-1}\log n$,
where $\kappa$ is a universal constant \underline{not} depending on $n$, $M$, $N$ or $x$.
Moreover, we also have the concentration inequality 
\[
\pp\{|F_H(x)-\ee(F_H(x))| \ge t\} \le 2\exp\biggl(-\frac{nt^2}{2\kappa\log n}\biggr)
\]
for every $t\ge 0$, where $\kappa$ is the same constant as in the variance bound.
\end{thm}
A remarkable aspect of Theorem \ref{free} is that the constant $\kappa$ is just a numerical constant independent of everything else. Note also that $H\mapsto F_H(x)$ is a discontinuous map. We believe that such a result cannot be established via gaussian type concentration of measure for orthogonal and unitary matrices (Gromov \& Milman \cite{gromovmilman83} and Szarek \cite{szarek90}). 

Incidentally, Voiculescu used the results of Gromov-Milman and Szarek in his celebrated work \cite{voiculescu91} that connected free probability theory with random matrices. That is an example of an area where concentration results such as Theorem \ref{free} may be relevant.

\subsection{The main result}
Let $G$ be a compact topological group. Then there exists a $G$-valued random variable $X$ with the properties that for any $x\in G$, the random variables $xX$, $Xx$ and $X^{-1}$ all have the same distribution as $X$. The distribution of $X$ is called the (normalized) Haar measure on $G$. The existence and uniqueness of the Haar measure is a classical result (see e.g.\ Rudin \cite{rudin73}, Theorem 5.14).
Let $Y$ be another $G$-valued random variable with the following properties:
\begin{enumerate}
\item The random variable $Y^{-1}$ has the same distribution as $Y$; that is, the law of $Y$ is {\it symmetric}.
\item For any $x\in G$, $xYx^{-1}$ has the same distribution as $Y$. In other words, the distribution of $Y$ is `constant on the conjugacy classes of~$G$'.
\end{enumerate} 
Recall that for two random variables $U$ and $V$ taking value in some separable space $\xx$, the supremum of $|\pp(U\in B) - \pp(V\in B)|$ as $B$ ranges over all Borel subsets of $\xx$ is called the total variation distance between the laws of $U$ and $V$, often denoted simply by~$d_{TV}(U,V)$. 
\begin{thm}\label{gpaction}
Let $G, X, Y$ be as above, with $X$ and $Y$ independent. Let $f:G \ra \rr$ be a bounded measurable function such that $\ee f(X) = 0$. Let $\|f\|_\infty = \sup_{x \in G} |f(x)|$ and 
\[
\|f\|_Y := \sup_{x\in G} \bigl[\ee(f(x)-f(Yx))^2\bigr]^{1/2}.
\]
Let $Y_1,Y_2,\ldots,$ be i.i.d.\ copies of $Y$. Suppose $a$ and $b$ are two positive constants such that $d_{TV}(Y_1Y_2\cdots Y_k, X) \le ae^{-bk}$ for every $k$, where $d_{TV}$ is the total variation metric. Let $A$ and $B$ be two numbers such that $\|f\|_\infty \le A$ and $\|f\|_Y \le B$. Let
\[
C = \frac{B^2}{b}\biggl[\biggl(\log \frac{4aA}{B}\biggr)^+ + \frac{b}{1-e^{-b}}\biggr].
\]
Then $\var(f(X))\le C/2$, and for any $t\ge 0$, $\pp\{|f(X)|\ge t\} \le 2e^{-t^2/C}$. 
\end{thm}
The main term in the bound is $B^2/b$; the term within the brackets will always contribute just a `factor of $\log n$' in applications (see discussion  in the next subsection).

Recall that if $ae^{-bt}$ expresses the correct rate of decay of the total variation distance, then $\tau := b^{-1}\log a$ is the mixing time of the Markov chain. Thus, the theorem roughly says the following: the deviation of $f(X)$ from its mean is of the order of $B\sqrt{\tau}$, where $B$ is a bound on the size of $f(x)-f(Yx)$, and $\tau$ is the mixing time of the Markov chain induced by $Y$. 

\subsection{Outline of the proofs}\label{outline}
Given a reversible Markov kernel $P$ and a function $f$ on the state space, the function
\begin{equation}\label{fdefn}
F(x,y) := \sum_{k=0}^\infty (P^k f(x) - P^k f(y))
\end{equation}
has the properties that $F(x,y) = -F(y,x)$, and 
\[
\ee(F(X_0,X_1)|X_0) = f(X_0) - \ee f(X_0),
\]
where $X_0, X_1,\ldots$ is a stationary Markov chain from the kernel $P$. Using these two properties and some intuition from Stein's method, we show that
\begin{equation}\label{varexpress}
\var (f(X_0)) = \frac{1}{2}\ee\bigl((f(X_0)-f(X_1)) F(X_0,X_1)\bigr).
\end{equation}
For a continuous Markov semigroup $(P_t)_{t\ge 0}$ with unique invariant measure $\mu$, the above identity is easily seen to be equivalent to 
\[
\var_\mu(f) = \int_0^\infty \mathcal{E}(f, P_t f) dt,
\]
where $\mathcal{E}$ is the Dirichlet form corresponding to the pair $((P_t)_{t\ge 0},\mu)$. Identities like this form the basis of the semigroup method for measure concentration.

Now suppose we can produce a number $B$ such that for all $k$, and all $x,y$ such that $y$ can be reached from $x$ in one step of the chain (i.e.\ $x$ and $y$ are `neighbors'), we have
\begin{equation}\label{naivebd}
|P^k f(x) - P^k f(y)| \le B.
\end{equation}
As $k$ increases beyond the mixing time $\tau$, $P^k f(x) - P^k f(y)$ vanishes exponentially fast. Combining, we see from the definition \eqref{fdefn} of $F$ that
\begin{equation}\label{naivebd2}
|F(x,y)|\lesssim \tau B.
\end{equation}
Using \eqref{naivebd2} and \eqref{naivebd} with $k=0$ in \eqref{varexpress}, we get
\[
\var(f(X_0)) \lesssim \frac{B^2\tau }{2}.
\]
This is the essence of the variance bound in Theorem \ref{gpaction}. The concentration inequality is obtained along a similar line. 

In practice, an inequality like \eqref{naivebd} is not very easy to establish. In fact, in the proof of Theorem \ref{gpaction} we are only able to prove \eqref{naivebd} in an average sense. The `constant on conjugacy classes' condition imposed on the random walk is required for our proof of \eqref{naivebd}. The key idea is to construct a coupling such that if two chains are started at neighboring sites, they continue to be on neighboring sites at each step. 

{\it The log factor.} The $\log$ factor in Theorem \ref{gpaction} arises from the $\log n$ terms appearing in the mixing times of random walks. In the above sketch, we used the fact that $P^k f(x) - P^k f(y)$ vanishes rapidly beyond $k > \tau$. Now, if $x$ and $y$ are neighboring states, this vanishing probably happens quicker, typically in $n$ steps instead of $n\log n$. But this is not stated in the standard theorems on Markov chain mixing. Any result in this direction (e.g.\ via path coupling) will suffice to remove the $\log$ factor from the statement of Theorem \ref{gpaction}.

The proof of Theorem \ref{free} is a direct application of Theorem \ref{gpaction}, proceeding as follows. First, we fix $x\in \rr$ and two matrices $M$ and $N$ as in the statement of Theorem \ref{free}. It is not difficult to see that the law of $F_H(x)$ is the same as that of $F_{UMU^* + N}(x)$, where $U$ is a Haar distributed unitary matrix. Accordingly, the state space is taken to be the set of $n\times n$ unitary matrices $\mathcal{U}_n$, and the function $f$ is defined as
\[
f(U) = F_{UMU^* + N}(x).
\]
We consider a random walk on $\mathcal{U}_n$ generated by conjugation with certain {\it random reflections}. The total variation rate of convergence to equilibrium for this walk is directly available from the literature \cite{porod96}.

\subsection{Discussion of existing literature} 
There is not much literature on the concentration of Haar measures. One early result is due to Maurey \cite{maurey79}, who investigated the Haar measure on the group $S_n$ of permutations of $n$ elements. The setting in Maurey's theorem is a particular case of ours, with $Y$ being a random transposition of two elements. 

Maurey's result was generalized in the lecture notes of Milman and Schechtman (\cite{milmanschechtman86}, Theorem 7.12) using a martingale argument.  Talagrand, in his famous treatment \cite{talagrand95}, made a substantial improvement on Maurey's result that allows one to go beyond `bounded differences'. The recent paper of Luczak and McDiarmid \cite{luczakmcdiarmid03} is also worthy of note.

The other group that has been studied for concentration of measure is the special orthogonal group $SO_n$, i.e., the group of $n\times n$ orthogonal matrices with determinant~$1$. 
The chief result about the concentration of Haar measure on this group is due to Gromov \& Milman~\cite{gromovmilman83}. As mentioned before, this result was used by Voiculescu \cite{voiculescu91} is his work connecting random matrix theory with free probability. 

However, other than the results about $S_n$ and $SO_n$ mentioned above, there is very little of general theory about the concentration of Haar measures. Theorem \ref{gpaction} is possibly the first result of its kind, and also the first result that connects rates of convergence to stationarity of random walks on groups with concentration of the invariant measures. Random walks on groups have received extensive attention following the pioneering works of Diaconis and Shahshahani \cite{diaconis81} and Diaconis and Saloff-Coste \cite{ds96}. Theorem \ref{gpaction} allows us to translate results about the rate of convergence to stationarity of random walks on groups which are `constant on conjugacy classes' to concentration inequalities under the Haar measure. Indeed, we will use one such available result \cite{porod96} to obtain the concentration of the Haar measure on the group of unitary matrices of order $n$ with respect to the rank distance for $n\times n $ matrices (the rank distance is defined as $d(A,B) := \mathrm{rank}(A-B)$).

Finally, let us clarify that the `concentration property of groups' as defined by Gromov \& Milman \cite{gromovmilman83} and investigated by Pestov (see, e.g.\ \cite{pestov02}) is not related to the sort of things that we are investigating.

\section{Proofs}\label{proof}
\subsection{Proof of Theorem \ref{gpaction}}
We begin with the observation that $Y$ defines a reversible Markov kernel $P$ in a natural way: For any $f:G\ra \rr$ such that $\ee|f(X)|<\infty$, let
\begin{equation}\label{ykernel}
Pf(x) := \ee f(Yx) = \ee f(xx^{-1}Yx) = \ee f(xY).
\end{equation}
The reversibility of this kernel can be proved as follows: Since $yX$ has the same distribution as $X$ for any $y\in G$, and $X,Y$ are independent, therefore $Y$ and $YX$ are also independent. Also, $Y^{-1}$ has the same distribution as $Y$. Hence, the pair $(X,Y)$ has the same distribution as $(YX, Y^{-1})$. Consequently, the pairs $(X,YX)$ and $(YX, Y^{-1}YX) = (YX,X)$ also have the same distribution. In other words, $(X,YX)$ is an exchangeable pair of random variables. This is equivalent to saying that $P$ is a reversible Markov kernel. The following lemma gives the most important information about this kernel that we require.
\begin{lmm}\label{lm1}
Under the hypothesis of Theorem \ref{gpaction}, and with $P$ defined in~\eqref{ykernel}, we have 
\begin{align*}
&\sum_{k=0}^\infty \ee|(f(x)-f(Yx)(P^kf(x)-P^kf(Yx))| \\
 &\le \frac{B^2}{b}\biggl[\biggl(\log \frac{4aA}{B}\biggr)^+ + \frac{b}{1-e^{-b}}\biggr].
\end{align*}
\end{lmm}
\begin{proof}
Note that for any $x\in G$,
\begin{align*}
|P^k f(x)| &= |P^k f(x) - \ee f(X)| = |\ee f(Y_1\cdots Y_k x) - \ee f(Xx)|\\
&\le 2\|f\|_\infty d_{TV}(Y_1\cdots Y_k, X) \le 2\|f\|_\infty a e^{-bk}.
\end{align*}
This shows, in particular, that for any $x\in G$, we have
\begin{equation}\label{sumcond}
\sum_{k=0}^\infty |P^k f(x)|\le \frac{2\|f\|_\infty a }{1-e^{-b}} < \infty.
\end{equation}
More importantly, it gives the bound
\begin{equation}\label{min1}
\begin{split}
&\ee|(f(x)-f(Yx))(P^kf(x)-P^kf(Yx))| \\
&\le 4\|f\|_\infty a e^{-bk} \ee|f(x)-f(Yx)| \le 4\|f\|_\infty a e^{-bk} \|f\|_Y.
\end{split}
\end{equation}
Now recall the assumption 2 that for any $y\in G$, $y^{-1}Yy$ has the same distribution as $Y$. Thus, for any $x,y\in G$,
\[
Pf(yx) = \ee f(Yyx) = \ee f(y y^{-1}Yyx) = \ee f(yYx).
\] 
So, if we let $Y^\prime$ be an independent copy of $Y$, then
\begin{align*}
\ee(Pf(x)-Pf(Yx))^2 &= \ee(\ee(f(Y^\prime x) - f(YY^\prime x)|Y)^2) \\
&\le \ee(f(Y^\prime x) - f(YY^\prime x))^2 \\
&\le \sup_{y^\prime\in G} \ee(f(y^\prime x) - f(Yy^\prime x))^2 = \|f\|_Y^2.
\end{align*}
Thus, $\|Pf\|_Y \le \|f\|_Y$. Continuing by induction, we get $\|P^k f\|_Y \le \|f\|_Y$. This gives
\begin{align}\label{min2}
&\ee|(f(x)-f(Yx)(P^kf(x)-P^kf(Yx))| \nonumber \\
&\le \bigl(\ee(f(x)-f(Yx))^2\bigr)^{1/2}\bigl(\ee(P^kf(x)-P^kf(Yx))^2\bigr)^{1/2} \nonumber \\
&\le \|f\|_Y \|P^k f\|_Y \le \|f\|_Y^2.
\end{align}
Using \eqref{min1} and \eqref{min2}, we get
\begin{align}\label{vx1}
&\ee|(f(x)-f(Yx)(P^kf(x)-P^kf(Yx))| \nonumber \\
&\le \sum_{k=0}^\infty \min\{\|f\|_Y^2, 4a \|f\|_\infty\|f\|_Y e^{-bk}\}\nonumber\\
&\le \sum_{k=0}^\infty \min\{B^2, 4a ABe^{-bk}\}\nonumber = B^2\sum_{k=0}^\infty \min\{1, 4aAB^{-1} e^{-bk}\}.
\end{align}
We shall now compute a bound on the above sum. For ease of notation let $\beta = 4aAB^{-1}$, and let 
$\gamma = b^{-1}\log \beta$.
If $\beta < 1$, the sum is just a geometric series which is easy to evaluate. Now assume $\beta\ge 1$. Then $\gamma$ is nonnegative. Now, an easy verification shows that $\beta e^{-b\gamma} = 1$, and $1 \ge \beta e^{-bk}$ if and only if $k \ge \gamma$. Let $k_0$ be the integer such that $k_0-1< \gamma \le k_0$. Then
\[
\sum_{k=0}^\infty \min\{1, \beta e^{-bk}\} \le k_0+ \sum_{k\ge k_0} \beta e^{-bk} = k_0 + \frac{\beta e^{-bk_0}}{1-e^{-b}}.
\]
Now the function
\[
g: x\mapsto x+\frac{\beta e^{-bx}}{1-e^{-b}}
\]
is convex and is therefore upper bounded by $\max\{g(\gamma), g(\gamma+1)\}$ on the interval $[\gamma,\gamma+1]$. A simple verification now shows that
\[
g(\gamma)=g(\gamma+1) = \gamma + \frac{1}{1-e^{-b}}.
\]
This completes the proof of the lemma. 
\end{proof}
\noindent Now define the function $F:G^2 \ra \rr$ as
\[
F(x_1,x_2) := \sum_{k=0}^\infty (P^k f(x_1)-P^k f(x_2)),
\]
where $f$ is the function under consideration in Theorem \ref{gpaction}.
By \eqref{sumcond}, the sum converges everywhere. The following lemma establishes the relevant properties of $F$.
\begin{lmm}\label{lm2}
The function $F$ satisfies $F(x_1,x_2)=-F(x_2,x_1)$ \[
\ee(F(x,Yx)) = f(x).
\]
\end{lmm}
\begin{proof}
The first property is obvious. Now, 
$\ee(P^k f(Yx)) = P^{k+1}f(x)$. 
Thus, for any $N$, we have
\[
\sum_{k=0}^N \ee(P^kf(x)-P^k f(Yx)) = f(x) - P^{N+1}f(x).
\] 
Now, by \eqref{sumcond}, we have $\lim_{N\ra \infty} P^{N+1} f(x)=0$. The uniform bound in \eqref{sumcond} also allows us to use the dominated convergence theorem. This completes the proof of the lemma.
\end{proof}
\noindent We are now ready to finish the proof of Theorem \ref{gpaction}. First, let us define the function
\begin{align*}
v(x) &:= \ee|(f(x)-f(Yx))F(x,Yx)|. 
\end{align*}
By Lemma \ref{lm2} and the independence of $X$ and $Y$, we get
\begin{equation}\label{step1}
\ee (f(X)^2) = \ee(f(X)F(X,YX)).
\end{equation}
Since $F(x_1,x_2)\equiv -F(x_2,x_1)$, we also get $\ee(f(X)^2)=-\ee(f(X)F(YX,X))$. Now, as proved in the beginning of this section, $(X,YX)$ is an exchangeable pair of random variables. Thus, 
\begin{equation}\label{step2}
\ee(f(X)^2) = - \ee(f(X)F(YX,X)) = -\ee(f(YX)F(X,YX)).
\end{equation}
Combining \eqref{step1} and \eqref{step2}, we get
\[
\ee(f(X)^2) =\frac{1}{2}\ee\bigl((f(X)-f(YX))F(X,YX)\bigr) \le \frac{1}{2} \ee( v(X)).
\]
By Lemma \ref{lm1}, $|v(x)|\le C$ for each $x$, where $C$ is as defined in the statement of Theorem~\ref{gpaction}. This proves the second moment bound. 

For the exponential inequality, let us define $\varphi(\theta) := \ee(e^{\theta f(X)})$ for each $\theta \in \rr$. Since $f$ is a bounded function, therefore $\varphi$ is differentiable and 
\[
\varphi^\prime (\theta) = \ee(e^{\theta f(X)} f(X)) = \ee(e^{\theta f(X)} F(X,YX)).
\]
Proceeding exactly as before, we get
\[
\varphi^\prime(\theta)=\frac{1}{2} \ee\bigl((e^{\theta f(X)} - e^{\theta f(YX)})F(X,YX)\bigr).
\]
Now, for any $u,v\in \rr$, we have
\begin{align*}
\biggl|\frac{e^u-e^v}{u-v}\biggr| &=\int_0^1 e^{tu + (1-t)v}dt \le \int_0^1 (te^u+(1-t)e^v)dt 
= \frac{1}{2}(e^u + e^v).
\end{align*}
Using this, and the exchangeability of $(X,XY)$ and the symmetry of $|F|$, we get
\begin{align*}
|\varphi^\prime(\theta)| &\le \frac{|\theta|}{4}\ee((e^{\theta f(X)} + e^{\theta f(YX)})|(f(X)-f(YX)) F(X,YX)|) \\
&\le \frac{|\theta|}{2} \ee(e^{\theta f(X)} |(f(X)-f(YX)) F(X,YX)|)\\
&= \frac{|\theta|}{2} \ee(e^{\theta f(X)} v(X)) \le \frac{C|\theta| \varphi(\theta)}{2}.
\end{align*}
This gives $\varphi(\theta)\le C\theta^2/4$ for all $\theta$. The proof can now be easily completed via routine arguments.

\subsection{Proof of Theorem \ref{free}}
Throughout this subsection, $\uu_n$ will denote the group of unitary matrices of order $n$.
To prove Theorem \ref{free}, we first need to establish a theorem about the concentration of the Haar measure on $\uu_n$. Existing results of the type discussed in Section \ref{intro} cannot give concentration bounds for $F_H$, since they are based on the Hilbert-Schmidt distance, which is not suitable for this purpose. Instead, we shall work with the {\it rank distance}, defined as $d(M,N) := \mathrm{rank}(M-N)$. The empirical distribution function is well-behaved with respect to this metric, as shown by the following lemma of Bai \cite{bai99}:
\begin{lmm}\label{bailmm}
\textup{[Bai \cite{bai99}, Lemma 2.2]}
Let $M$ and $N$ be two $n\times n$ hermitian matrices, with empirical distribution functions $F_M$ and $F_N$. Then 
\[
\|F_M-F_N\|_\infty \le \frac{1}{n}\mathrm{rank}(M-N).
\]
\end{lmm}
This lemma is an easy consequence of the interlacing inequalities for eigenvalues of hermitian matrices. (It seems possible that this already existed in the literature before \cite{bai99}, but we could not find any reference.)

To find the concentration of the Haar measure on $\uu_n$ with respect to the rank distance, we need a random walk which takes `small steps' with respect to this metric. 

Let $G = \uu_n$ and $X$ be a Haar-distributed random variable on $\uu_n$. We define the r.v.\ $Y$ required for generating the random walk for Theorem \ref{gpaction} as follows: Let $Y= I-(1-e^{i\varphi})uu^*$, where $u$ is drawn uniformly from the unit sphere in $\cc^n$, and $\varphi$ is drawn independently from the distribution on $[0, 2\pi)$ with density proportional to $(\sin(\varphi/2))^{n-1}$. Multiplication by $Y$ represents a random reflection across a randomly chosen subspace.
It is easy to verify that $Y\in \uu_n$. Now, $Y^{-1}= Y^* = I - (1-e^{-i\varphi})uu^* = I-(1-e^{i(2\pi - \varphi)})uu^*$ has the same distribution as $Y$, since $2\pi -\varphi$ has the same distribution as $\varphi$. Also, for any $U\in \uu_n$,
\[
UYU^* = I- (1-e^{i\varphi})(Uu)(Uu)^*,
\]
and $Uu$ is again uniformly distributed over the unit sphere in $\cc^n$. Hence $Y$ satisfies all the properties required for Theorem \ref{gpaction}. 

Following a sketch of Diaconis \& Shahshahani \cite{diaconis86}, Ursula Porod \cite{porod96} proved the following result about the rate of convergence to stationarity of the random walk induced by $Y$:
\begin{thm}
\textup{[Porod \cite{porod96}]}
Let $X,Y$ be as above. Let $Y_1,Y_2,\ldots,$ be i.i.d.\ copies of $Y$. There exists universal constants $\alpha, \beta, c_0$, such that whenever $n \ge 16$ and $k \ge \frac{1}{2}n\log n + c_0n$, we have
\begin{equation}\label{alphabeta}
d_{TV}(Y_1\cdots Y_k, X) \le \alpha n^{\beta/2} e^{-\beta k/n},
\end{equation}
where $d_{TV}$ denotes the total variation distance.
\end{thm}
Substituting $k=\frac{1}{2}n\log n + c_0 n$, we get $\alpha e^{-\beta c_0}$ on the right hand side. Thus by suitably increasing $\alpha$ such that $\alpha e^{-\beta c_0} \ge 1$, we can drop the condition that $k\ge \frac{1}{2}n\log n + c_0n$. 
Combining Porod's theorem with Theorem~\ref{gpaction}, we get the following result about concentration of the Haar measure on $\uu_n$.
\vskip.2in
\begin{prop}\label{unitary}
Let $G=\uu_n$ and $X,Y$ be as above, with $n\ge 16$. Let $f:\uu_n \ra\rr$ be a function such that $\ee f(X)=0$. Let $\|f\|_Y = \sup_{U\in \uu_n} [\ee(f(U)-f(YU))^2]^{1/2}$. Let $A$ and $B$ be constants such that $\|f\|_\infty \le A$ and $\|f\|_Y \le B$. Let 
\[
C = \frac{nB^2}{\beta}\biggl[\biggl(\log \frac{4\alpha n^\beta A}{B}\biggr)^+ + \frac{\beta/n}{1-e^{-\beta/n}}\biggr],
\]
where $\alpha$ and $\beta$ are as in \eqref{alphabeta}.
Then $\var(f(X))\le C/2$, and for any $t\ge 0$, we have $\pp\{|f(X)|\ge t\} \le 2e^{-t^2/C}$.
\end{prop}
\vskip.3in
We are now ready to finish the proof of Theorem \ref{free}.
\begin{proof}[Proof of Theorem \ref{free}]
Suppose $U$ and $V$ are independent Haar-distributed random unitary matrices of order $n$, and $H=UMU^* + VNV^*$. 
The matrix $V^*HV = V^*UM U^*V + N$ has the same spectrum as $H$. Also, $V^*U$ is again Haar distributed. Hence, we can consider, without loss of generality, the spectrum of 
\[
H = XM X^* + N,
\]
where $X$ follows the Haar distribution on $\uu_n$.  Now let
\[
H^\prime = (YX) M (YX)^* + N.
\]
Recall that $Y = I-(1-e^{i\varphi})uu^*$, where $u$ is drawn from the uniform distribution on the unit sphere in $\cc^n$, and $\varphi$ is drawn independently from the distribution on $[0,2\pi)$ with density proportional to $(\sin(\varphi/2))^{n-1}$. Let $\delta = 1-e^{i\varphi}$. Then 
\begin{align*}
H-H^\prime &= XM X^* - (I-\delta uu^*)XM X^*(I-\bar{\delta} uu^*)\\
&=\delta Huu^* + \bar{\delta} uu^*H - |\delta|^2 uu^*Huu^*.
\end{align*}
The three summands are all of rank $1$ and thus $\mathrm{rank}(H-H^\prime) \le 3$. Thus by Lemma \ref{bailmm}, we see that
\begin{equation}\label{fdiff}
\|F_H - F_{H^\prime}\|_\infty\le \frac{3}{n}. 
\end{equation}
Now fix a point $x\in \rr$, and let $f:\uu_n \ra \rr$ be the map which takes $X$ to $F_H(x)$. Then by~\eqref{fdiff}, we have
\[
|f(X)-f(YX)| \le \frac{3}{n} \ \ \text{ for all possible values of $X$ and $Y$.}
\]
Thus, $\|f\|_Y \le 3/n$. Also, $\|f\|_\infty \le 1$. Thus, in Proposition \ref{unitary}, we get $C \le \kappa \log n + c$ for some universal constants $\kappa$ and $c$. By choosing $\kappa$ large enough, we can drop the assumption that $n\ge 16$ and also put $c=0$. 
This completes the proof. 
\end{proof}
\vskip.3in
\noindent{\bf Acknowledgments.} The author thanks the referee for several useful comments and corrections, and a small improvement in Theorem \ref{gpaction}. Sincere thanks are also due to Persi Diaconis, Yuval Peres, and Assaf Naor for reading the manuscript and providing advice.

\begin{small}

\end{small}

\end{document}